\documentclass[11pt]{article}
\usepackage{amsmath,amsthm,amsfonts,amssymb,mathrsfs,epsfig}
\newtheorem{theorem}{Theorem}
\newtheorem{lemma}{Lemma}

\newtheorem{proposition}{Proposition}
\theoremstyle{remark}

\def\N{\mathbb{N}}
\def\Z{\mathbb{Z}}
\def\R{\mathbb{R}}

\def\PP{\mathsf{P}}

\def\AA{\mathscr{A}}

\def\SS{\mathscr{S}}

\def\t{\mathfrak{t}}

\renewcommand{\phi}{\varphi}
\renewcommand{\epsilon}{\varepsilon} 
\newcommand{\1}{{\text{\Large $\mathfrak 1$}}}
\newcommand{\comp}{\raisebox{0.1ex}{\scriptsize $\circ$}}

\renewcommand{\theta}{\vartheta}

\newcommand{\keywords}[1]{ \noindent {\footnotesize
             {\small \em Keywords and phrases.} {\sc #1} } }
\newcommand{\ams}[2]{  \noindent {\footnotesize
             {\small \em AMS {\rm 2000} subject classifications.
             {\rm Primary {\sc #1}; secondary {\sc #2}} } } }

\begin{document}

\title{\bf Stationary flows 
and uniqueness of invariant measures}
\author{
\begin{minipage}{5.5cm}
\begin{center}
{\sc Fran\c{c}ois Baccelli} \\
\'Ecole Normale Sup\'erieure
\\
Paris, France
\end{center}
\end{minipage}
\hfill
\begin{minipage}{5.5cm}
\begin{center}
{\sc Takis Konstantopoulos} \\
Heriot-Watt University\\
Edinburgh, UK
\end{center}
\end{minipage}
}
\date{}
\maketitle

\begin{abstract}
In this short paper, we consider a quadruple $(\Omega, \AA, \theta, \mu)$,
where $\AA$ is a $\sigma$-algebra of subsets of $\Omega$,
and $\theta$ is a measurable bijection from $\Omega$ into itself
that preserves the measure $\mu$.
For each $B \in \AA$, we consider the 
measure $\mu_B$ obtained by taking cycles (excursions)
of iterates of $\theta$ from $B$. 
We then derive a relation for $\mu_B$
that involves the forward and backward hitting times of $B$
by the trajectory $(\theta^n \omega, n \in \Z)$ at
a point $\omega \in \Omega$.
Although classical in appearance, its use in obtaining uniqueness
of invariant measures of various stochastic models seems to be new.
We apply the concept to countable Markov chains and Harris
processes. 

\vspace*{2mm}

\vspace*{2mm}
\ams{28D05,60G10}{60J10,60J05}

\end{abstract}

\section{Introduction}
This paper was initiated from the following question.
It is classical that,
for a Markov chain $(X_n, n \ge 0)$ with a countable state space $S$ that 
possesses a positive recurrent state $b \in S$, 
there is at least one invariant probability measure
$\pi^{(b)}$ on $S$ which is defined by the usual ``cycle formula'':
\[
\pi^{(b)}(A) = \frac{1}{E_b \t_b} E_b \sum_{n=0}^{\t_b-1} \1(X_n \in A),
\]
where $\t_b$ is the first return time to $b$.
To ensure that $\pi^{(b)}$ is the only invariant probability measure
we need, in addition, to ensure that the only positive recurrent states
are those that communicate with $b$ 
(this holds, for instance, if the chain is irreducible). 
There are several proofs of uniqueness, ranging from analytic 
(e.g.\ by means of the Perron-Frobenius
theorem which itself can be proved in a number of ways--see, e.g.\
Lind and Marcus (1995) for a geometric proof) to probabilistic
(e.g.\ by means of applying the Doeblin coupling construction: this requires,
in addition, aperiodicity--see, e.g.\ Thorisson (2000)).
The question we posed is whether there is a way to prove uniqueness directly
from the way that $\pi^{(b)}$ is constructed by the cycle formula.
If so, can we do this for Markov chains in a general state space?
And finally, how ``Markovian'' is the proof of uniqueness (can the
``local'' character of definition of $\pi_b$ be extended to
other processes)?

In answering the question, we abstracted the problem and lifted
it to a general measurable space $(\Omega, \AA)$ endowed with
a measurable bijective transformation $\theta$ that preserves
some measure $\mu$.
The point of view appears to be new, although
the tools used below are quite natural in Ergodic Theory and in the
construction of Palm Probabilities. 
The origin of these tools can be traced, as far as we can tell, to a
paper by Kac (1947).
In Section \ref{mastersection} we define, for each $B \in \AA$, the 
forwards and backwards hitting times of $B$ by the iterates of $\theta$
(called $T_B, \widetilde T_B$, respectively) and the measure
\[
\mu_B(A) = \int_B d \mu ~ \sum_{n=0}^{T_B-1} \1_{\theta^{-n} A}.
\]
Theorem \ref{mftheorem} states the basic formula of interest:
\[
\mu_B(A) = \mu(A, \widetilde T_B < \infty).
\]
It can be read as: on the event that $B$ has been visited in the past
at least once, the measures $\mu_B$ and $\mu$ coincide.
Thus, if $\mu(B) > 0$, Poincar\'e's recurrence lemma 
(recalled as Lemma \ref{poinc}), $\mu_B=\mu$ for all $B$.
In Section \ref{mcsection}, we consider a Markov chain on
a countable set $S$. Assuming irreducibility and positive recurrence,
the previous observation immediately yields a unique probability
measure $\pi$ on $S$ such that $\pi \PP = \pi$, which answers the original
question. Finally,
in Section \ref{harrissection}, we consider a Harris chain
and show uniqueness of the invariant probability measure constructed by means of cycles away from
a recurrent regeneration set $R$.

\section{The master formula}
\label{mastersection}

Let $(\Omega, \AA)$ be a measurable space
and $\theta: \Omega \to \Omega$ a measurable bijection.
For $A, B \in \AA$ define the following functions:
\begin{subequations}
\label{basic}
\begin{align}
\label{basic1}
T_B & \equiv T_B(\omega; \theta) := \inf\{n \in \N:~ \theta^n\omega \in B\}\\
\label{basic2}
M_B(A) & \equiv M_B(A, \omega; \theta) 
:= \sum_{0 \le n < T_B(\omega; \theta)} \1(\theta^n\omega \in A),
\end{align}
\end{subequations}
stressing that both take values in $\N \cup \{+\infty\}:=
\{1,2,\ldots\} \cup \{+\infty\}$,
where $\inf \varnothing = +\infty$.
The definition of $M_B(A)$ requires giving a meaning to
the quantity $\theta^{T_B}$. We let
$\Omega_B = \{T_B < \infty\}$, and 
define $\theta^{T_B} : \Omega_B \to \Omega_B$ by
\[
\big(\theta^{T_B}\big)(\omega) := \theta^{T_B(\omega)}(\omega), \quad
\omega \in \Omega_B.
\]
On $\Omega-\Omega_B$, we define $\theta^{T_B}$ rather arbitrarily,
e.g.\ by letting it to be the identity on it.
We can easily see that $\theta^{T_B}$ is invertible
with $\big(\theta^{T_B}\big)^{-1} = \theta^{-T_B}$, 
and where $\theta^{-T_B}$ is defined in a similar way.
We shall also need \eqref{basic1}-\eqref{basic2}
when using $\theta^{-1}$ in place of 
$\theta$:
\begin{align*}
\widetilde T_B & \equiv
T_B(\omega; \theta^{-1}) := \inf\{n \in \N:~ \theta^{-n} \omega \in B\}\\
\widetilde M_B(A) & \equiv M_B(A, \omega; \theta^{-1}) := \sum_{0 \le n < T_B(\omega; \theta^{-1})} \1(\theta^{-n} \omega \in A).
\end{align*}
The interpretation is that $M_B(A)$ evaluated at $\omega$ 
is the number of times the
forward trajectory $(\omega, \theta\omega, \theta^2 \omega,\ldots)$
visits the set $A$ up to (and not including) the time it visits
the set $B$. Similarly, $\widetilde M_B(A)$ refers to the backward
trajectory $(\omega, \theta^{-1}\omega, \theta^{-2} \omega,\ldots)$.
There is a certain ``duality'' between $M_B(A)$ and $\widetilde T_B$ on
one hand and $\widetilde M_B(A)$ and $T_B$ on the other, once
we integrate against an invariant measure. We discuss this next.
Recall first the following standard lemma:
\begin{lemma}[Poincar\'e recurrence]
\label{poinc}
If the measure $\mu$ on $(\Omega, \AA)$
is preserved by $\theta$ then, for all $B \in \AA$,
\begin{equation}
\label{pr}
\mu(B) = \mu(B, \widetilde T_B < \infty)
= \mu(B, T_B < \infty).
\end{equation}
\end{lemma}
\proof
This follows from
\begin{multline*}
\mu(B^c, \widetilde T_B=\infty) = \lim_{n\to\infty}
 \mu(B^c \cap \theta B^c \cap \cdots \cap \theta^{n-1} B^c)
\\
= \lim_{n\to\infty}  \mu(\theta B^c \cap \theta^2 B^c \cap \cdots \cap \theta^n B^c)
=\mu(\widetilde T_B=\infty),
\end{multline*}
and similarly for $T_B$.
\qed

\noindent
In other words, $\widetilde T_B<\infty$ and
$T_B<\infty$, $\mu$-a.e.\ on $B$. This is used in proving:
\begin{theorem}
\label{mftheorem}
If the measure $\mu$ on $(\Omega, \AA)$
is preserved by $\theta$, then, for all $A, B \in \AA$,
\begin{subequations}
\label{mf}
\begin{align}
\mu_B(A) &:= \int_B M_B(A)  d\mu = \int_A \1(\widetilde T_B<\infty)  d\mu,
\label{mf1}\\
\widetilde \mu_B(A) 
&:= \int_B \widetilde M_B(A) d\mu = \int_A \1(T_B<\infty)  d\mu.
\label{mf2}
\end{align}
\end{subequations}
\end{theorem}
\proof
We only need to show the first identity.
\begin{align}
\int_B M_B(A)  d\mu 
&= 
\mu(A \cap B) +
\sum_{n \ge 1} \mu(B \cap \theta^{-1}B^c \cap \cdots \cap \theta^{-n} B^c
\cap \theta^{-n} A)
\nonumber \\
&= \mu(A \cap B) +
\sum_{n \ge 1} \mu(\theta^{n} B \cap \theta^{n-1}B^c \cap \cdots \cap B^c \cap A)
\nonumber \\
&= \mu(A \cap B) + \sum_{n \ge 1} \mu(\widetilde T_B=n, A \cap B^c)
\nonumber \\
&= \mu(A \cap B) + \mu(A \cap B^c, \widetilde T_B < \infty)=
\mu(A, \widetilde T_B < \infty),
\label{arg}
\end{align}
where the Poincar\'e recurrence formula (and more precisely its
consequence that $\mu(A\cap B, \widetilde T_B = \infty)=0$)
was used to obtain the last equality.
\qed

\begin{proposition}[strong invariance]\label{prop.invar}
If the measure $\mu$ on $(\Omega, \AA)$
is preserved by $\theta$, then, 
its restriction $\mu(\cdot \cap B)$ on some $B \in \AA$ is preserved
by $\theta^{T_B}$ and by $\theta^{\widetilde T_B}$, i.e.,
for all $A, B \in \AA$,
\[
\mu(B \cap \theta^{-T_B} A) 
= \mu(B \cap \theta^{-\widetilde T_B} A) 
= \mu(A \cap B).
\]
\end{proposition}
\noindent
{\em Note:} The terminology {\em strong invariance} is by analogy 
to the strong Markov property.
\proof
Since, due to the Poincar\'e recurrence, 
$\mu(B \cap \theta^{-T_B} A \cap (\Omega-\Omega_B))=0$, we have
\begin{align*}
\mu(B \cap \theta^{-T_B} A)
&= \sum_{n=1}^\infty \mu(B \cap \theta^{-n}A, T_B=n) \\
&= \sum_{n=1}^\infty \mu(B \cap \theta^{-n}A \cap \theta^{-1}  B^c \cap
\cdots \cap \theta^{-(n-1)} B^c \cap \theta^{-n} B) \\
&=  \sum_{n=1}^\infty \mu(B  \cap \theta^{-1}  B^c \cap
\cdots \cap \theta^{-(n-1)} B^c \cap \theta^{-n} (A\cap B)) \\
&=  \sum_{n=1}^\infty \mu(\theta^n B \cap \theta^{n-1} B^c \cap
\cdots \cap \theta B^c \cap A \cap B) \\
&= \sum_{n=1}^\infty \mu(\widetilde T_B=n, A \cap B) = \mu(\widetilde T_B < \infty,  A \cap  B)
=\mu(A \cap B),
\end{align*}
where the latter equality again follows from the Poincar\'e 
recurrence \eqref{pr}.
The second assertion is proved in the same manner.
\qed
\begin{proposition}
\label{strongprop}
If the measure $\mu$ on $(\Omega, \AA)$
is preserved $\theta$, then, for all $B \in \AA$,
the measures $\mu_B(\cdot)$, $\widetilde \mu_B(\cdot)$, 
defined by \eqref{mf1}, \eqref{mf2},
respectively, are also preserved by $\theta$.
\end{proposition}
\proof 
Note that
\[
M_B(\theta^{-1} A) - M_B(A) = \1_{\theta^{-T_B}A} - \1_{A}.
\]
Using this and Proposition \ref{prop.invar} we obtain
\[
\mu_B(\theta^{-1} A) - \mu_B(A)
= \int_B M_B(\theta^{-1} A)  d\mu - \int_B M_B(A)  d\mu
= \mu(B \cap \theta^{-T_B}A) - \mu(B \cap A) = 0.
\]
\qed

\paragraph{Some remarks:}

\paragraph{(i)}
Since $M_B(\Omega) = T_B$,
$\widetilde M_B(\Omega) = \widetilde T_B$,
we have, from Theorem \ref{mftheorem},
\begin{equation}
\label{qqq}
\int_B T_B  d\mu = \mu(\widetilde T_B < \infty),
\quad
\int_B \widetilde T_B  d\mu = \mu(T_B < \infty).
\end{equation}
Thus, if $\mu=P$ is a probability measure and if $E$ denotes
integration with respect to $P$,
then
\[
E T_B \1_B = P(\widetilde T_B < \infty) \le 1.
\]
If, in addition, $P(B) > 0$ then $P(\widetilde T_B < \infty)
\ge P(\widetilde T_B < \infty, B) = P(B) > 0$
and so
\[
E(T_B \mid B) = \frac{1}{P(B \mid \widetilde T_B < \infty)},
\]
where, as usual, $E(T_B \mid B) = \frac{E T_B \1_B }{P(B)}$.
This is slightly more general than Kac' formula (Kac (1947)).
Similar formula holds, of course, for $E(\widetilde T_B \mid B)$:
\[
E(\widetilde T_B \mid B) = \frac{1}{P(B \mid T_B < \infty)}.
\]

\paragraph{(ii)}
Since $T_B \ge 1$, $\widetilde T_B \ge 1$, \eqref{qqq} implies that:
\[
\mu(B) \le \mu(\widetilde T_B < \infty) \wedge \mu(T_B < \infty).
\]
This leads to the following equivalences:
\[
\mu(B) > 0 
\iff \mu(T_B < \infty) >0
\iff \mu(\widetilde T_B < \infty) >0.
\]
Indeed, if $\mu(B) > 0$ then $\mu(\widetilde T_B < \infty) >0$
from the last inequality.
Conversely, if $\mu(B) =0$ then \eqref{arg} shows that
$\mu(\widetilde T_B < \infty) =0$.

\paragraph{(iii)}
The function $B \mapsto \mu_B(A)$ can be thought of as
a pre-capacity. Indeed, let 
\[
\Psi(\omega) := \{\theta^{-1} \omega, \theta^{-2} \omega, \ldots\}
\]
and consider it as a random set.
Then
\[
\mu_B(A) = P(A, \Psi \cap B \not = \varnothing)
\]
is the pre-capacity functional of the random set $\Psi$ 
(see Molchanov (2005).)
We avoid using the terminology capacity because there no topological
properties of $\Psi$ are introduced.
An interesting problem would be to investigate properties
of the function $\mu_B(A)$ jointly in $A, B$.

\paragraph{(iv)}
Theorem \ref{mftheorem} and Proposition \ref{strongprop}
should of course be linked to the cycle formula of Palm calculus,
and Proposition \ref{prop.invar} to the invariance of the Palm measure.
The main point here is that within this discrete time setting,
there is no need to invoke the general theory (Baccelli and Br\'emaud (2003)).

\paragraph{(v)}
Some results do not require the invertibility of $\theta$.
For instance, $\mu(B) = \mu(B, T_B < \infty)=
\mu(B, \cup_{n=1}^{\infty} \theta^{-n} B)$ holds for any $\mu$-preserving
measurable map $\theta$ (Lemma \ref{poinc}).
However, the main formulae \eqref{mf1}-\eqref{mf2} that exhibit
the ``duality'' between forward and backward iterates of $\theta$,
do require invertibility.
On the other hand, even without using Theorem 
\ref{mftheorem} and Propositions \ref{prop.invar}-\ref{strongprop}, 
we can show that the measure
$\nu_B(A):=\mu(A, T_B<\infty)$ satisfies $\nu_B(\theta A) = \nu_B(A)$
directly.
To do this, note that $\{T_B \comp \theta < \infty\} = B \cup \{T_B < \infty\}$
and write
\begin{align*}
\mu(\theta A,  T_B < \infty)
&= \mu(A,  T_B \comp \theta < \infty)\\
&= \mu((A\setminus B) \cup (A \cap B),~ B \cup \{T_B < \infty\})\\
&= \mu(A\setminus B, T_B < \infty) + \mu(A\cap B, T_B < \infty) \\
&=  \mu(A, T_B < \infty).
\end{align*}
This, incidentally, gives a second proof of Proposition \ref{strongprop}.

\section{Uniqueness in Markov chains}
\label{mcsection}
Suppose that $\PP=[p_{i,j}]$ is a stochastic matrix on a  countable 
state space $S = \{a,b,c,\ldots,i,j,\ldots\}$, i.e.
\[
p_{i,j} \ge 0,
\quad
\sum_{k \in S} p_{i,k}=1,
\quad
i,j \in S.
\]
Assume that it is
\begin{itemize}
\item[(i)]
irreducible (each $i$ {\em communicates} with each $j$ in $S$),
\item[(ii)]
positive recurrent (starting from some $i$ the expected return time
to $i$ has finite expectation).
\end{itemize}
These properties depend entirely on the matrix $\PP$.
It is classical that:
\begin{theorem}
\label{class}
If (i) and (ii) hold then
there is a unique probability $\pi$ on $S$ such that
$\pi \PP = \pi$.
\end{theorem}
We wish to show this by using the idea developed in the previous section.

\paragraph{\em Proof of existence}
It is uniqueness that is novel here.
Existence of such a $\pi$ is immediately answered by the ``cycle
formula'': Let $(X_0, X_1, \ldots)$ be a realisation of the 
Markov chain with transition probability matrix $\PP$.
Fix some state $b$, let
\[
\t_b := \inf\{n \ge 1: X_n=b\},
\]
and define the probability $\pi^{(b)}$ on $S$ by
\begin{equation*}
\pi^{(b)}(a) = \frac{E_b  \sum_{n=0}^{\t_b-1} \1(X_n=a) }
{E_b \t_b},
\quad a \in S,
\end{equation*}
where $E_b$ is expectation conditional on $X_0=b$.
That this $\pi^{(b)}$ is an invariant probability measure (satisfies
$\pi^{(b)} \PP = \pi^{(b)}$ is standard (see, e.g.\ Br\'emaud (1999)). 
It is important to note that $\pi^{(b)}$ depends entirely on
the stochastic matrix $\PP$ only.
\qed

\paragraph{\em Proof of uniqueness}
To show uniqueness, we work at the level of sequences, i.e.\
with the space $\Omega = S^\Z$, whose elements are
denoted by $\omega=(\omega_n, n \in \Z)$, equipped with the cylinder
$\sigma$-algebra $\AA$. 
We consider the natural shift
\[
\theta: (n \mapsto \omega_n) \mapsto (n \mapsto \omega_{n+1}),
\]
which is obviously $\AA$-measurable and invertible. We are thus
in the setup of the earlier section.
Consider a probability $\pi$ on $S$ satisfying $\pi \PP = \pi$,
and let $P$ be the probability measure on $(S^\Z, \AA)$
defined by
\begin{multline}
P(\{\omega\in \Omega:~\omega_m=i_m, \ldots, \omega_n=i_n\}) 
= \pi(i_m) p_{i_m,i_{m+1}} \cdots p_{i_{n-1},i_n},\\
\quad i_m, \ldots, i_n \in S, 
\quad m, n \in \Z, \quad m \le n.
\end{multline}
Consider also the random variables 
\[
X_n(\omega) := \omega_n, \quad
\omega \in \Omega, \quad
n \in \Z.
\]
Under $P$, the sequence $(X_n)$ is a Markov
chain with transition probability matrix $\sf P$.
Clearly, the measure $P$ is preserved by $\theta$
and, by Proposition \ref{strongprop}, 
so are the measures 
\[
P_B(\cdot) = E \1_B M_B(\cdot) = \int_B M_B(\cdot) dP,
\]
where $T_B, M_B$ are given by \eqref{basic1}-\eqref{basic2},
for any $B \in \AA$.
Fix some $b \in S$, and consider the set
\[
B= \{\omega\in \Omega:~\omega_0=b\}.
\]
Observe that
\[
\t_b(\omega) := \inf\{n \ge 1: \omega_n=b\} = T_B(\omega),
\quad
\widetilde \t_b(\omega) := \inf\{n \ge 1: \omega_{-n}=b\} = \widetilde T_B(\omega).
\]
By Theorem \ref{mftheorem},
\begin{equation}    \label{here}
P_B(A) = E \1_B M_B(A) = P(A, \widetilde T_B<\infty)
=P(A, \widetilde \t_b < \infty),
\quad A \in \AA.
\end{equation}
By (i) and (ii) we have $P(\t_b < \infty)=1$,
$P(\widetilde \t_b < \infty)=1$,
and so \eqref{here} yields
\begin{equation*}    
P_B(A) = P(A), \quad A \in \AA,
\end{equation*}
and $E\1_B T_B = 1$.
Therefore,
\[
P(A) = P_B(A)
= E \1_B M_B(A) = \frac{E \1_B M_B(A)}{E \1_BT_B} =
\frac{E_b \sum_{n=0}^{\t_b-1} \1(\theta^n \omega \in A)}{E_b \t_b}.
\]
So, if we pick 
\[
A:=\{\omega\in \Omega:~ \omega_0=a\},
\]
we conclude that $\pi(a) = \pi^{(b)}(a)$ for all $a \in S$.
Thus, an arbitrary invariant probability measure $\pi$ must be equal to 
the specific measure $\pi^{(b)}$; whence the uniqueness.
\qed

\paragraph{Remarks:}
\paragraph{(i)}
The last argument directly proves that
\[
\frac{E_b  \sum_{n=0}^{\t_b-1} \1(X_n=a) }{E_b \t_b}
=
\frac{E_c  \sum_{n=0}^{\t_c-1} \1(X_n=a) }{E_c \t_c},
\]
the so-called exchange formula of (discrete-index) Palm theory
(see also Konstantopoulos and Zazanis (1995)).
\paragraph{(ii)}
Only the existence proof used the Markov property. The
uniqueness proof was at the level of stationary processes.
\paragraph{(iii)} 
In essence, uniqueness follows from the following two facts:
\begin{itemize}
\item Unique determination of the Palm measure: thanks to the Markov
setting considered here, the Palm law of a cycle starting from a
given state until the chain returns to this state is uniquely
determined by the transition matrix; 
\item Slivnyak's inverse construction: this construction shows that the
stationary law of a point process is fully determined by its Palm measure.
(See Slivnyak (1962).)
\end{itemize}
Again, the main point here is that there is no need to invoke
the general theory.
\paragraph{(iv)}
The same argument can be used to show the weaker result:
\begin{theorem}
\label{class2}
Suppose that (ii) holds (every state is positive recurrent)
Let $S=\cup_{i \ge 1} S_i$, be the
decomposition of $S$ into its irreducible components.
Let $b_i \in S_i$, for all $i \ge 1$.
Then every probability $\pi$ on $S$ such that
$\pi \PP = \pi$ is a convex combination of the measures
$\pi^{(b_i)}$.
\end{theorem}

\section{Uniqueness in Harris chains}
\label{harrissection}
The method explained above can also be applied to yield
a proof of uniqueness for the invariant probability measure
of a positive Harris recurrent chain.

A Markov process $(X_n)$ with values in a Polish space $(S, \SS)$
and transition kernel 
\[
K(x, \cdot) = P_x(X_1 \in \cdot)
\]
is called
{\em Harris recurrent} or, simply, Harris chain
(Asmussen (2003)) if it possesses a recurrent regeneration set $R \in S$.
This means that
\begin{itemize}
\item[(i)]
\[
P_x(\t_R < \infty) = 1, \quad x \in S,
\]
where
\[
\t_R := \inf\{n \in \N: X_n \in R\};
\]
\item[(ii)]
there is a probability measure $\lambda$ on $(S,\SS)$, an $\epsilon > 0$,
and $\ell \in \N$, such that
\[
K^\ell(x, \cdot) \ge \epsilon \lambda(\cdot), \quad x \in \R,
\]
where 
\[
K^\ell(x,\cdot) = P_x(X_\ell \in \cdot).
\]
\end{itemize}
The chain is called {\em positive Harris recurrent} if, in addition
to (i) and (ii) we also have
\begin{itemize}
\item[(iii)]
\[
E_\lambda \t_R < \infty,
\]
\end{itemize}
where, as usual, $E_\lambda$ denotes expectation with respect
to $P_\lambda(\cdot) := \int_S \lambda(dx) P_x(\cdot)$.

We here give a proof of the following:
\begin{theorem}
\label{harrisuniq}
A positive Harris recurrent chain possesses a unique invariant probability
measure.
\end{theorem}

Note that this theorem is proved in the paper of Athreya and Ney (1978)
by different methods and only in the case $\ell=1$. There is a substantial
difference between the $\ell=1$ and $\ell>1$ cases in that the cycles
defined by the iterates of the stopping time $\t$ (see \eqref{ct} below)
are not independent.

{\em Proof of Theorem \ref{harrisuniq}.}
Existence is standard (see Asmussen (2003))
and requires construction of the chain on 
a suitable probability space. We repeat the construction here.
In addition to the chain, we consider a sequence $(\zeta_n)$ of
i.i.d.\ Bernoulli random variables taking values $1$ or $0$ with
probability $\epsilon$ or $1-\epsilon$ respectively.
Informally, whenever $X_n \in R$ distribute
$X_{n+\ell}$ according to $\lambda$ if $\zeta_n=1$
or according to 
$\displaystyle \frac{K^\ell(x,\cdot)-\epsilon \lambda(\cdot)}{1-\epsilon}$
if $\zeta_n=0$, and, conditional on $(X_n, X_{n+\ell})$, distribute
$(X_{n+1}, \ldots, X_{n+\ell-1})$ by respecting the given Markov kernel.
Otherwise, if $X_n \not \in R$, then ignore $\zeta_n$ and continue
the chain as usual.
Formally, we define an $\ell$-th order Markov chain $(X_n, \zeta_n)$
with values in $S \times \{0,1\}$ via the following:
Let $G(dx_1, \ldots, dx_{\ell-1} | x, y)$ be the conditional
distribution of $(X_{1}, \ldots, X_{\ell-1})$ given that $X_0=x, X_\ell=y$, i.e.
\begin{equation}
\label{conddist}
G(dx_1, \ldots, dx_{\ell-1} | x, y)
:= 
\frac{K(x,dx_1) \cdots 
K(x_{\ell-2}, d x_{\ell-1}) K(x_{\ell-1}, d y)}
{\int_{S^{\ell-1}} K(x,dx_1') \cdots
K(x_{\ell-2}', d x_{\ell-1}') K(x_{\ell-1}', d y)},
\end{equation}
where the integration in the denominator is with respect to the variables
$x_1', \ldots, x_{\ell-1}'$ and the ratio is
to be understood as a Radon-Nikod\'ym derivative with respect
to $y$. 
Then let
\begin{multline}
P(X_{n+i} \in dx_i, 1 \le i\le\ell \mid X_n=x, \zeta_n=\alpha)
\\
=
\begin{cases}
\lambda(d x_\ell)
G(dx_1, \ldots, dx_{\ell-1} | x, x_\ell),
& \text{ if } x \in R, ~ \alpha =1 \\
\displaystyle{\frac{K^\ell(x,dx_\ell)-\epsilon \lambda(dx_\ell)}{1-\epsilon}}~
G(dx_1, \ldots, dx_{\ell-1} | x, x_\ell),
& \text{ if } x \in R, ~ \alpha =0 
\\
K(x,dx_1) \cdots K(x_{\ell-1}, dx_\ell),
& \text{ otherwise}
\end{cases}
\label{both}
\end{multline}
and finally require that, for all $n$,
\begin{multline}
P(X_{n+i} \in dx_i, \zeta_i=\alpha_i, 1 \le i\le\ell\mid X_m, \zeta_m, m \le n)
\\
=
p(\alpha_1) \cdots p(\alpha_\ell) ~
P(X_{n+i} \in dx_i, 1 \le i\le\ell \mid X_n, \zeta_n),
\label{mp}
\end{multline}
where $\alpha_i \in \{0,1\}$ and $p(0) := \epsilon$, $p(1):=1-\epsilon$.

It is easy to see that $(X_n)$ is a realisation of the Harris chain
with the given transition kernel $K$, and that $(\zeta_n)$ is
an i.i.d.\ sequence; the two sequences are dependent.

Consider
\begin{equation}
\label{ct}
\t:= \inf\{n : X_{n-\ell} \in R,~ \zeta_{n-\ell}=1\},
\end{equation}
(so that $X_{\t}$ has distribution $\lambda$)
and define
\begin{equation}
\label{fromlambda}
\pi(\cdot) :=
\frac{E_\lambda \sum_{n=0}^{\t-1} \1(X_n \in \cdot)}
{E_\lambda \t},
\end{equation}
It is now standard to check that $\pi(\cdot)$
is an invariant probability measure for the chain $(X_n)$.

To prove uniqueness, we shall again consider the same construction
defined by \eqref{conddist},
\eqref{both} and \eqref{mp}, and, in addition, we shall assume that
the chain is stationary and therefore defined over the index set $\Z$.
Specifically, our probability space is $\Omega = (S\times \{0,1\})^\Z$,
equipped with the natural cylinder $\sigma$-algebra $\AA$. A typical
element of $\Omega$ is denoted by $\omega = \big((x_n, \zeta_n), n \in \Z\big)$.
The shift is again the natural one:
\[
\theta: (n \mapsto (x_n, \zeta_n)) \mapsto (n \mapsto (x_{n+1}, \zeta_{n+1})).
\]
The probability measure $P$ on $(\Omega, \AA)$ is such that
it makes the coordinate process an $\ell$-th order
Markov chain with transition kernel 
defined through \eqref{conddist}, \eqref{both}
and \eqref{mp}, and is invariant under $\theta$. 
(Thus, we have created a setup $(\Omega, \AA, \theta, P)$,
as in Section \ref{mastersection}, where $P$ plays the r\^ole of $\mu$
and, here, $P(\Omega)=1$.)
We now prove that there can be only one such $P$. To this end, let
\[
B:= \{\omega=(x,\zeta) \in \Omega:~ x_{-\ell}\in R,~ \zeta_{-\ell}=1\}.
\]
By our assumptions, $P(T_B < \infty)=1$, $P(\widetilde T_B < \infty)=1$.
By Theorem \ref{mftheorem},
\[
P_B(A) = E \1_B M_B(A) = P(A, \widetilde T_B< \infty) = P(A), \quad
A \in \AA,
\]
and $E\1_B T_B =1$. 
But 
\[
P_B(A) = \frac{E \1_B M_B(A)}{E\1_B T_B}
= \frac{\displaystyle E\left[\sum_{n=0}^{\t-1} 
\1(\theta^n \omega \in A) \mid 
x_{-\ell} \in R, \zeta_{-\ell}=1 \right]}
{\displaystyle E[\t \mid x_{-\ell} \in R, \zeta_{-\ell}=1 ]}
=
\frac{\displaystyle E_\lambda \sum_{n=0}^{\t-1}
\1(\theta^n \omega \in A)}{E_\lambda \t},
\]
since, by construction, $P(C | x_{-\ell} \in R, \zeta_{-\ell}=1 )
= P_\lambda(C)$ for any $C$ in the $\sigma$-algebra
generated by $(\omega_n, n \ge 0)$.
Taking $A:=\{\omega=(x,\zeta)\in \Omega:~ x_0 \in \cdot\}$
we conclude that any $\theta$-invariant probability measure $P$
that preserves the given Markovian structure must have a
marginal given by \eqref{fromlambda}.
This proves uniqueness.
\qed

Final note: The proof of uniqueness, again, uses arguments
that do not rely on the Markov property.
As such, it would be worth exploiting it further in
stochastic scenaria with absence of Markovian property.

\section*{References}
\begin{description}
\item
{\sc Asmussen, S.} (2003).
{\em Applied Probability and Queues},
2nd ed.
Springer-Verlag. 
\item
{\sc Athreya, K.B.\ and Ney, P.} (1978).
A new approach to the limit theory of recurrent Markov chains.
{\em Trans.\ Amer.\ Math.\ Soc.} {\bf 245}, 493-501.
\item
{\sc Baccelli, F.\ and Br\'emaud, P.} (2003).
{\em Elements of Queueing Theory.}
Springer-Verlag.
\item
{\sc Br\'emaud, P.} (1999).
{\em Markov Chains: Gibbs Fields, Monte Carlo Simulation,
and Queues.}
Springer Verlag.
\item
{\sc Kac, M.} (1947).
On the notion of recurrence in discrete stochastic processes.
{\em Bull.\ AMS} {\bf 53}, 1002-1010.
\item
{\sc Konstantopoulos, T. and Zazanis, M.}  (1995). 
A discrete time proof of Neveu's exchange formula.
{\em J.\ Appl.\ Probability} {\bf 32}, 917-921.
\item
{\sc Lind, D. and Marcus, B.} (2000).
{\em An Introduction to Symbolic Dynamics and Coding.}
Cambridge University Press.
\item
{\sc Molchanov, I.} (2005).
{\em Theory of Random Sets.}
Springer-Verlag.
\item
{\sc Slivnyak, I.M.} (1962).
Some properties of stationary flows of homogeneous random events.
{\em Th.\ Prob.\ Appl.} {\bf 7}, 336-341.
\item
{\sc Thorisson, H.} (2000).
{\em Coupling, Stationarity, and Regeneration}.
Springer.
\end{description}

\vspace*{1cm}
\noindent
\hfill
\begin{minipage}[t]{6cm}
\small \sc
Authors' addresses:\\[3mm]
Fran\c{c}ois Baccelli\\
D\'epartement d'Informatique\\
\'Ecole Normale Sup\'erieure\\
45 rue d'Ulm\\
F-75230 Paris Cedex 05, France\\
E-mail: {\tt Francois.Baccelli@ens.fr}
\\[5mm]
Takis Konstantopoulos\\
School of Mathematical Sciences\\
Heriot-Watt University\\
Edinburgh EH14 4AS, UK\\
E-mail: {\tt takis@ma.hw.ac.uk}
\end{minipage}

\end{document}